\documentclass[graybox]{svmult}

\usepackage{type1cm}        % activate if the above 3 fonts are
\usepackage{makeidx}         % allows index generation
\usepackage{graphicx}        % standard LaTeX graphics tool
\usepackage{multicol}        % used for the two-column index
\usepackage[bottom]{footmisc}% places footnotes at page bottom
\usepackage{multirow}

\usepackage{newtxtext}
\usepackage{newtxmath}       % selects Times Roman as basic font

\usepackage{color}
\usepackage{ulem}
\usepackage{pifont}

\usepackage{pgfplots}
\usepackage{tikz}
\usetikzlibrary{arrows.meta,calc,shapes,fit,positioning,shadows,backgrounds}
\tikzset{no shadows/.style={general shadow/.style=}}
\def\myscale{0.5}
\tikzset{block/.style={rectangle split, draw, rectangle split parts=2,
text width=14em*\myscale, text centered, rounded corners, minimum height=4em*\myscale},
grnblock/.style={materia, fill=green!20, text width=10em*\myscale, text centered, rounded corners, minimum height=4em*\myscale}, 
 materia/.style={draw, fill=white!20, text width=15.0em*\myscale, text centered, minimum height=1.5em*\myscale,drop shadow},
whtblock/.style={no shadows, materia, text width=30em*\myscale, minimum width=20em*\myscale, minimum height=10em*\myscale, rounded corners, drop shadow},    
line/.style={draw, -{Latex[length=2mm*\myscale,width=1mm]}},
cloud/.style={draw, ellipse,fill=white!20, node distance=3cm*\myscale,    minimum height=4em*\myscale},  
container/.style={draw, rectangle,dashed,inner sep=0.9cm*\myscale, rounded
corners,fill=orange!10,minimum height=2cm*\myscale},
container2/.style={draw, rectangle,dashed,inner sep=0.9cm*\myscale, rounded
corners,fill=blue!10,minimum height=2cm*\myscale},
container3/.style={draw, rectangle,dashed,inner sep=0.9cm*\myscale, rounded
corners,fill=green!10,minimum height=2cm*\myscale}, 
box/.style={draw, rectangle, rounded corners, thick, node 
distance=7em*\myscale, 
text width=6em*\myscale, text centered, minimum height=3.5em*\myscale}, 
every node/.style={font=\scriptsize}}
\tikzstyle{arrow} = [thick,->,>=stealth]

\usepackage{mathtools}
\usepackage{listings}
\makeatletter
\DeclareRobustCommand\ttfamily
        {\not@math@alphabet\ttfamily\mathtt
         \fontfamily\ttdefault\selectfont\hyphenchar\font=-1\relax}
\makeatother
\DeclareTextFontCommand{\mytexttt}{\ttfamily\hyphenchar\font=45\relax}

\def\nonlinschwarzop{\mytexttt{NonLinearSchwarzOperator}}
\def\coarsenonlinschwarzop{\mytexttt{CoarseNonLinearSchwarzOperator}}
\def\schwarzop{\mytexttt{SchwarzOperator}}
\def\simplecoarseop{\mytexttt{SimpleCoarseOperator}}
\def\coarseop{\mytexttt{CoarseOperator}}
\def\overlappingop{\mytexttt{OverlappingOperator}}
\def\combineop{\mytexttt{CombineOperator}}
\def\nonlincombineop{\mytexttt{NonLinearCombineOperator}}
\makeindex             % used for the subject index
\usepackage[capitalise]{cleveref}

\begin{document}

\title*{Nonlinear Two-Level Schwarz Methods:\\A Parallel Implementation in FROSch}
\author{Alexander Heinlein\orcidID{0000-0003-1578-8104}\\ Kyrill Ho\orcidID{0009-0001-4522-5647}\\ Axel Klawonn\orcidID{0000-0003-4765-7387} \\ Martin Lanser\orcidID{0000-0002-4232-9395}}
\institute{Alexander Heinlein$^{1}$, Kyrill Ho$^{2}$, Axel Klawonn$^{2,3}$, Martin Lanser$^{2,3}$ \at $^1$Delft University of Technology, Delft Institute of Applied Mathematics, Mekelweg 4, 2628 CD Delft, Netherlands e-mail: a.heinlein@tudelft.nl, \at $^2$Department of Mathematics and Computer Science, Division of Mathematics, University of Cologne, Weyertal 86-90, 50931 Cologne, Germany, \email{axel.klawonn@uni-koeln.de, martin.lanser@uni-koeln.de}, {url:~\url{https://www.numerik.uni-koeln.de}} \at $^3$Center for Data and Simulation Science, University of Cologne, Germany, {url:~\url{https://www.cds.uni-koeln.de}}}
\maketitle

\abstract{
  Owing to the ability of nonlinear domain decomposition methods to improve the nonlinear convergence behavior of Newton's method, they have experienced a rise in popularity recently in the context of problems for which Newton's method converges slowly or not at all. This article introduces a novel parallel implementation of a two-level nonlinear Schwarz solver based on the FROSch (Fast and Robust Overlapping Schwarz) solver framework, part of Sandia's Trilinos library. First, an introduction to the key concepts underlying two-level nonlinear Schwarz methods is given, including a brief overview of the coarse space used to build the second level. Next, the parallel implementation is discussed, followed by preliminary parallel results for a scalar nonlinear diffusion problem and a 2D nonlinear plane-stress Neo-Hooke elasticity problem with large deformations.
}

\section{Introduction and Nonlinear Problems}
\label{sec:intro}
There is a growing body of work on nonlinear domain decomposition methods (DDMs). The interest in this class of methods can be attributed to their favorable properties for solving highly nonlinear problems, together with their inherent parallelizability. One class of nonlinear DDMs is the class of nonlinear Schwarz methods~\cite{schwarz1,schwarz5,raspen,mspin1,mspin2,schwarz6} on which we focus here.
Previous parallel implementations of nonlinear Schwarz methods have not taken advantage of the GDSW (Generalized Dryja-Smith-Widlund) type coarse spaces \cite{gdsw1}, for example \cite{parnlschwarz1} and \cite{parnlschwarz2}. The two-level implementations that incorporate this family of coarse space have been serial prototype codes implemented in high-level scripting languages such as MATLAB. Following the favorable results seen in \cite{schwarz1} for the two-level nonlinear Schwarz method with energy-minimizing coarse spaces, we have implemented a parallel two-level nonlinear Schwarz solver based on the FROSch (Fast and Robust Overlapping Schwarz) solver framework. FROSch is part of Trilinos and implements linear multilevel Schwarz methods with GDSW type coarse spaces. It has been shown to scale efficiently for linear and nonlinear problems as a preconditioner for an iterative linear solver \cite{frosch2}. In this article, we introduce the new nonlinear Schwarz implementation and discuss the first scalability results for a scalar nonlinear diffusion problem and a 2D nonlinear Neo-Hooke elasticity problem.

\paragraph{\bf Nonlinear Problems}
As stated above, the first of two nonlinear problems that we consider is a scalar nonlinear diffusion problem
\begin{align}
    \label{eq:nonlindiff}
    \begin{split}
        -\nabla((u^2+1)\nabla u) = 1 \; &\mathrm{in}\; \Omega,\\
        u = 0 \; &\mathrm{on} \; \partial \Omega,
    \end{split}
\end{align}
where $\Omega$ is the unit square. 
As a second model problem we consider a vector valued 2D nonlinear elasticity problem using a compressible plane-stress Neo-Hooke material model:
find $u\in (H_0^1(\Omega))^2$, such that
\begin{align}
    \label{eq:nonlinelas}
    \begin{split}
    -\mathrm{div}(P(F)) = f_{vol}\;&\mathrm{in}\;\Omega,\\
        u = g_D \;&\mathrm{on}\;\partial\Omega_D,\\
        n\cdot P(F) = g_N\;&\mathrm{on}\;\partial\Omega_N,
    \end{split}
\end{align}
where $F$ is the deformation gradient, and $\Omega_D$ and $\Omega_N$ represent the Dirichlet and Neumann boundaries, respectively. The first Piola-Kirchhoff stress tensor is
\begin{equation*}
  P(F) = \frac{E}{1(1+\nu)}(F-F^{-T}) + \frac{E\nu}{(1+\nu)(1-2\nu)}\mathrm{ln}(\mathrm{det}(F)F^{-T}).
\end{equation*}
We set the Poisson ratio to $\nu = 0.3$ and Young's modulus to $E = 210\,\mathrm{GPa}$. To model a simple beam, we choose the domain $\Omega$ to be a rectangle with $width=5\textrm{m}$ and $height=1\textrm{m}$. We fix the beam on either end by setting $g_D = 0$ and $g_N = 0$. We apply a volume force in the vertical direction by setting $f_{vol} = (0, -f_x)$ for some scalar $f_x > 0$.

\section{A Parallel implementation of Nonlinear Schwarz Methods}
\label{sec:nl}
In this section, we provide an overview of the expressions that describe the one-level nonlinear Schwarz method that was first introduced in \cite{schwarz5} and its extension to two-level methods described in \cite{schwarz1} and \cite{schwarz4}. In addition, we introduce our initial parallel implementation. We use the same definition for the subdomains $\Omega_i$, their corresponding function spaces $V_i$, the prolongation $P_i:V_i\mapsto V$ and restriction $R_i:V\mapsto V_i$ operators, the coarse space $V_0$ and associated operators $P_0:V_0\mapsto V$ and $R_0:V\mapsto V_0$ as used in \cite{schwarz1}. In our tests we use the RGDSW (Reduced GDSW) coarse space with the option $2.2$ as introduced in \cite{rgdsw} but with a slight modification; the values of the coarse basis functions on interface edges that end at a Dirichlet boundary are determined with the same inverse Euclidean distance formula used for internal interface edges. The original construction imposes a constant value of one on these interface edges.

To describe the one-level nonlinear Schwarz method we begin with an abstract discretized nonlinear problem 
\begin{equation}
    F(u) = 0
\end{equation}
and define local nonlinear correction terms $T_i(u) : V\mapsto V_i$ as the solution of the local nonlinear problems 
\begin{equation}
    \label{eq:localcorrection}
    R_iF(u-P_iT_i(u)) = 0,\; i = 1,\dots, N.
\end{equation}
Combining these nonlinear correction terms results in the alternative problem
\begin{equation}
    \label{eq:nonlinprec}
    \mathcal{F}_1(u) = \sum_{i = 1}^{N}P_iT_i(u) = 0.
\end{equation}

The one-level nonlinear Schwarz method results by solving equation \eqref{eq:nonlinprec}  using Newton's method. It is denoted ASPIN (Additive Schwarz Preconditioned Inexact Newton) in \cite{schwarz5} because the tangent $D\mathcal{F}_1(u)$ is evaluated approximately. Here we always evaluate the tangent exactly yielding the ASPEN (Additive Schwarz Preconditioned Exact Newton) method. We thus have two nested Newton's methods. In the \textit{outer loop} the alternative problem $\mathcal{F}_1(u) = 0$ is solved. In each iteration of the outer loop, the local nonlinear corrections $T_i(u^{(k)})$ must be computed to evaluate $\mathcal{F}_1(u^{(k)})$. This is done in \textit{local inner loops} using Newton's method on each overlapping subdomain $\Omega_i$ to solve the local nonlinear problems \eqref{eq:localcorrection}. 
The outer loop requires the evaluation of the tangent $D\mathcal{F}_1(u)$ in each iteration. Taking the derivative of equation \eqref{eq:localcorrection} results in the following expression for the tangent
\begin{equation}
    \label{eq:tangent}
    D\mathcal{F}_1(u) = \sum_{i=1}^{N}P_iDT_i(u) = \sum_{i=1}^{N}P_i(R_iDF(u_i)P_i)^{-1}R_iDF(u_i)
\end{equation}
with $u_i \coloneqq u-P_iT_i(u)$. To simplify the expression for the tangent we introduce 
\begin{equation}
    \label{eq:q}
    Q_i(u) \coloneqq P_i(R_iDF(u)P_i)^{-1}R_iDF(u)
\end{equation}
and write $D\mathcal{F}_1(u) = \sum_{i=1}^NQ_i(u_i)$. In this article, we prefer to work with a restricted Schwarz approach called RASPEN (Restricted ASPEN) which was described in \cite{schwarz7}. This method introduces restricted prolongation operators $\widetilde{P}_i:V_i\mapsto V,\; i = 1,\dots,N$, that are chosen to fulfill the partition of unity property 
\begin{equation}
    \label{eq:pou}
    \sum_{i = 1}^{N}\widetilde{P}_iR_i = I.
\end{equation}

To extend the one-level method, we define a coarse nonlinear correction $T_0(u): V\mapsto V_0$ as the solution of the coarse nonlinear problem 
\begin{equation}
    R_0F(u-P_0T_0(u)) = 0.
\end{equation}
Coupling of the nonlinear coarse correction with the local nonlinear corrections can be done in a variety of ways. Here we only describe the two possibilities for which we later show results. For further possibilities, see, e.g., \cite{schwarz1}. The first way of coupling is a simple additive approach in which the alternative problem is defined as 
\begin{equation}
    \mathcal{F}_a(u) = \sum_{i=1}^{N}{P}_iT_i(u) + P_0T_0(u)= 0.
\end{equation}
The second approach is a multiplicative hybrid approach, denoted $h1$ in \cite{schwarz1} but which we denote here simply as $h$ with
\begin{equation}
    \mathcal{F}_h(u) = \sum_{i=1}^{N}{P}_iT_i(u - P_0T_0(u)) + P_0T_0(u)= 0.
\end{equation}
The tangents of both methods are given by
\begin{align}
    D\mathcal{F}_a(u) = \sum_{i = 0}^{N}Q_i(u_i) \quad {\rm and} \quad D\mathcal{F}_h(u) = &\sum_{i=1}^{N}Q_i(v_i)(I-Q_0(u_0)) + Q_0(u_0)
\end{align}
with the previous definition of $u_i$ expanded to include $i = 0$ and $v_i \coloneqq u_0 - P_iT_i(u_0)$ for $i = 1,\dots, N$.

Our current parallel implementation is based on extending the functionality of the linear Schwarz operators of the Trilinos package FROSch. We refer to \cite{frosch1} and \cite{frosch2} for an overview of the various FROSch operators and an explanation of how GDSW coarse basis functions are built. Constructing RGDSW coarse basis functions uses the same mechanisms. Here, we give a brief overview of the current state of our implementation.

The nonlinear Schwarz solver requires finite element assembly routines to reassemble the tangent and residual periodically. Currently, we use the C++ package FEDDLib (Finite Element and Domain Decomposition Library) \cite{feddlib} for this and other preprocessing tasks such as mesh generation and distribution. The outer loop of the nonlinear Schwarz solver is simply an implementation of Newton's method in which the tangent system
\begin{equation}
    \label{eq:tangentsys}
    D\mathcal{F}_x(u^{(k)})\delta u = \mathcal{F}_x(u^{(k)})
\end{equation}
is solved in each iteration with $x\in\{1,a,h\}$. This requires evaluation of the nonlinear operator $\mathcal{F}_x(u)$ and formation of the tangent $D\mathcal{F}_x(u)$ at an arbitrary linearization point $u$.

In our implementation, we have separated the evaluation of the residual into various operators that all inherit from existing FROSch Schwarz operators. The one-level residual $\mathcal{F}_1$ is evaluated by the class \nonlinschwarzop{} which inherits from \schwarzop{}. The coarse component $P_0T_0(u)$ of the two-level residual is evaluated by the \coarsenonlinschwarzop{} class. The tangent is similarly compartmentalized into the \overlappingop{} class which is an original FROSch class, and the \simplecoarseop{} class which inherits from \coarseop{}.

To construct arbitrary two-level nonlinear Schwarz methods we have implemented a pair of abstract operators \nonlincombineop{} and \combineop{} that define the residual $\mathcal{F}_x(u)$ and the tangent $D\mathcal{F}_x(u)$, respectively. Specializations of these classes define the order in which local and coarse operators are combined. Employing this interface, we have implemented the three methods (one-level, two-level additive, and two-level hybrid) and can easily implement further operators for other multiplicative and hybrid two-level methods.

Using our current nonlinear Schwarz implementation is very similar to using a classical Newton-Krylov-Schwarz solver. The user must provide the assembly routines for the residual and tangent of the problem they wish to solve, configure the solver parameters, and start the solver. For now, the assembly routines must be implemented within a FEDDLib \mytexttt{NonLinearProblem} class object. This is because the assembly routines for the residual and tangent are currently expected to use the interface prescribed therein. We plan to make our implementation more accessible in the future by providing the same assembly routine interface as the Trilinos nonlinear solver package \mytexttt{NOX}.

The key difference to using a Newton-Krylov-Schwarz solver is that an overlapping mesh partitioning is required. The overlapping partition must be element based rather than node based, that is, use the dual graph of the mesh. This facilitates assembly on the overlapping subdomains since finite element assembly routines generally perform quadrature on a per-element basis. Since the linear Schwarz preconditioner works algebraically, it only requires a node-based overlapping partitioning. We have implemented routines in the FEDDLib to build an overlapping distributed dual graph using METIS/ParMETIS \cite{parmetis} and distribute the mesh correspondingly.

\section{Parallel Numerical Results}
The parallel results shown in this section were generated using the NHR supercomputer Fritz; see the Acknowledgements for further details.

To solve the global tangent system \eqref{eq:tangentsys} we use GMRES provided by the Trilinos package Belos \cite{belos}. We set the relative GMRES tolerance to $10^{-6}$, the maximum number of iterations to $1000$, and the maximum number of restarts to $20$. In the inner Newton methods, we use Intel MKL Pardiso as sparse direct solver \cite{pardiso}.

As a proof of concept to ensure that our implementation of the two-level nonlinear Schwarz variants is working as expected, we first test the simple nonlinear diffusion problem introduced in section \ref{sec:intro}. We do not compare the nonlinear Schwarz variants with a classical Newton-Krylov-Schwarz approach for this problem as it is too simple and we do not expect performance improvements from nonlinear Schwarz methods. We test weak parallel scalability of the one-level, the two-level additive, and the two-level hybrid nonlinear Schwarz methods. We use subdomains with $300$ elements in each dimension on a regular square mesh, set the relative tolerance of the outer Newton method to $10^{-4}$, the relative tolerance of the inner Newton methods to $10^{-5}$ and the absolute tolerance of the inner Newton methods to $10^{-11}$. In the two-level variants, we use the RGDSW coarse space \cite{rgdsw}. As an initial value, we set the solution to zero everywhere.

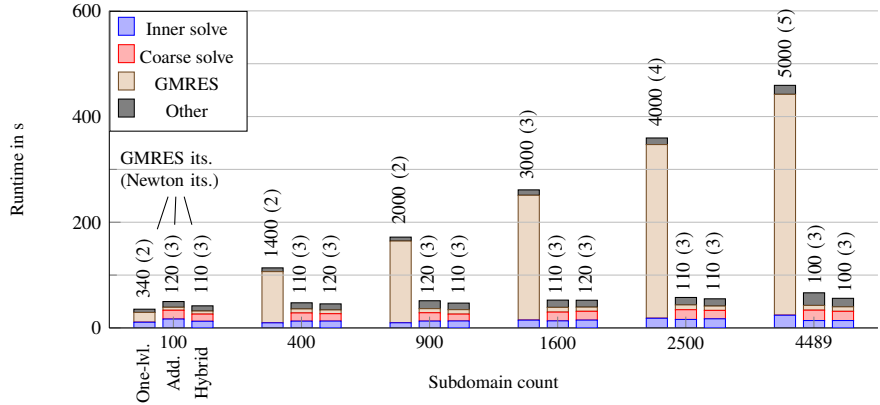
\begin{figure}
    \centering
    \begin{tikzpicture}[
    every axis/.style={
    legend style={at={(0.0,1)},anchor=north west},
    axis lines*=left, ymajorgrids, yminorgrids,
    width=11.8cm, height=5.8cm,
    ymin=0,
    ymax=600,
    ybar stacked,
    bar width=8pt,
    minor y tick num=1,
    xtick={1,2,3,4,5,6},
    xticklabels from table={\datatableonelevel}{SubdomainCount},
    xticklabel style={rotate=0,xshift=0ex,anchor=north},
    ylabel={Runtime in s},
    xlabel={Subdomain count},
    }
]
\pgfplotstableread{
Location SubdomainCount  GlobalSolve   InnerSolve   CoarseSolve   GMRES    Other
1        100             49.9          16.8         16.6          5.6      10.9
2        400             47.5          12.8         15.5          7.6      11.6
3        900             51.4          12.9         15.8          7.75     14.95
4        1600            52.6          13.3         16.6          8.9      13.8
5        2500            57.6          15.96        18.45         9.1      14.09
6        4489            66.4          13.95        19.6          8.95     23.9
}\datatableh
\pgfplotstableread{
Location SubdomainCount  GlobalSolve   InnerSolve   CoarseSolve   GMRES    Other
1        100             41.8          12.6         13.55         5.8      9.85
2        400             45.5          12.94        14.0          7.1      11.46
3        900             47.0          13.15        13.0          8.4      12.45
4        1600            52.4          14.8         16.65         7.9      13.05
5        2500            55.1          17.26        15.9          8.2      13.74
6        4489            56.02         14           17.4          8.5      16.12
}\datatableadd
\pgfplotstableread{
Location SubdomainCount  GlobalSolve   InnerSolve   CoarseSolve   GMRES    Other
1        100             35.4          10.74        0             18.5     6.16
2        400             113.6         9.6          0             96.9     7.1
3        900             171.8         9.7          0             154.7    7.4
4        1600            261.2         14.7         0             236.3    10.2
5        2500            359.4         18.24        0             328.7    12.46
6        4489            459           24           0             418      17
}\datatableonelevel

\begin{axis}[bar shift=-11pt, hide axis]
    \addplot+ table [x=Location, y=InnerSolve] {\datatableonelevel};
    \addplot+ table [x=Location, y=CoarseSolve] {\datatableonelevel};
    \addplot+ table [x=Location, y=GMRES] {\datatableonelevel};
    \addplot+ table [x=Location, y=Other] {\datatableonelevel};
\end{axis}

\begin{axis}[bar shift=0pt, hide axis]
    \addplot+ table [x=Location, y=InnerSolve] {\datatableh};
    \addplot+ table [x=Location, y=CoarseSolve] {\datatableh};
    \addplot+ table [x=Location, y=GMRES] {\datatableh};
    \addplot+ table [x=Location, y=Other] {\datatableh};
\end{axis}

\begin{axis}[bar shift=11pt]
    \addplot+ table [y=InnerSolve] {\datatableadd}; \addlegendentry{Inner solve}
    \addplot+ table [y=CoarseSolve] {\datatableadd}; \addlegendentry{Coarse solve}
    \addplot+ table [y=GMRES] {\datatableadd}; \addlegendentry{GMRES}
    \addplot+ table [y=Other] {\datatableadd}; \addlegendentry{Other}
\end{axis}

\node[rotate=90] at (0.45,-0.6) {One-lvl.};
\node[rotate=90] at (0.85,-0.6) {Add.};
\node[rotate=90] at (1.25,-0.6) {Hybrid};

\node[rotate=0, text width=1.5cm] (gmres) at (0.9,2.1) {GMRES its. (Newton its.)};

\node[rotate=90] (one) at (.47,0.80) {$340$ $(2)$};
\node[rotate=90] (two) at (0.85,0.86) {$120$ $(3)$};
\node[rotate=90] (three) at (1.24,0.86) {$110$ $(3)$};

\node[rotate=90] at (2.17,1.4) {$1400$ $(2)$};
\node[rotate=90] at (2.55,.86) {$110$ $(3)$};
\node[rotate=90] at (2.95,.86) {$120$ $(3)$};

\node[rotate=90] at (3.87,1.85) {$2000$ $(2)$};
\node[rotate=90] at (4.25,0.86) {$120$ $(3)$};
\node[rotate=90] at (4.66,0.86) {$110$ $(3)$};

\node[rotate=90] at (5.57,2.45) {$3000$ $(3)$};
\node[rotate=90] at (5.95,0.86) {$110$ $(3)$};
\node[rotate=90] at (6.35,0.86) {$120$ $(3)$};

\node[rotate=90] at (7.28,3.15) {$4000$ $(4)$};
\node[rotate=90] at (7.65,0.9) {$110$ $(3)$};
\node[rotate=90] at (8.05,0.9) {$110$ $(3)$};

\node[rotate=90] at (9.,3.8) {$5000$ $(5)$};
\node[rotate=90] at (9.4,0.99) {$100$ $(3)$};
\node[rotate=90] at (9.8,0.89) {$100$ $(3)$};

\draw [thin] (gmres) --  (one);
\draw [thin] (gmres) --  (two);
\draw [thin] (gmres) --  (three);

\end{tikzpicture}
    \caption{Time to solution of the one-level and two-level additive and hybrid nonlinear Schwarz variants for the nonlinear diffusion problem \eqref{eq:nonlindiff}. The outer Newton iterations and total number of GMRES iterations summed over all outer iterations is shown above each column.}
    \label{fig:ttslaplace}
\end{figure}

Figure \ref{fig:ttslaplace} demonstrates the necessity of a coarse level even for a simple nonlinear problem. The number of outer nonlinear iterations is not affected by whether a coarse space is used or not and remains fairly constant around three iterations. However, the condition number of the tangent $D\mathcal{F}(u)$ deteriorates with more subdomains in the one-level case and results in an increasing number of GMRES iterations. The increase in the time to solution of the one-level method clearly results directly from the increase in the number of GMRES iterations. The perfect weak scalability of the two-level variants is probably due to the fact that the nonlinear diffusion problem is simple enough that the coarse correction generates a very good approximation of the solution. We expect to see an increase in the time to solution of the two-level variants when the coarse problem grows large enough.

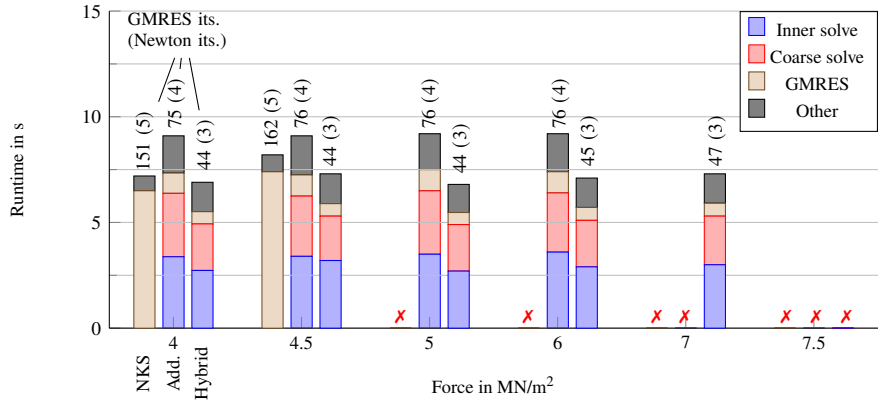
\begin{figure}
    \centering
    \begin{tikzpicture}[
    every axis/.style={
    legend style={at={(1,1)},anchor=north east},
    axis lines*=left, ymajorgrids, yminorgrids,
    width=11.8cm, height=5.8cm,
    ymin=0,
    ymax=15,
    ybar stacked,
    bar width=8pt,
    minor y tick num=1,
    xtick={1,2,3,4,5,6},
    xticklabels from table={\hybrid}{Force},
    xticklabel style={rotate=0,xshift=0ex,anchor=north},
    ylabel={Runtime in s},
    xlabel={Force in MN/m$^{2}$},
    }
]

\pgfplotstableread{
Location Force  GlobalSolve   InnerSolve   CoarseSolve   GMRES    Other
1        4      6.9           2.73         2.2           0.58     1.39
2        4.5    7.3           3.2          2.1           0.59     1.41
3        5      6.8           2.7          2.2           0.57     1.33
4        6      7.1           2.9          2.2           0.61     1.39
5        7      7.3           3            2.3           0.61     1.39
6        7.5    0             0            0             0        0
}\hybrid

\pgfplotstableread{
Location Force  GlobalSolve   InnerSolve   CoarseSolve   GMRES    Other
1        4      9.1           3.38         3             0.96     1.76
2        4.5    9.1           3.4          2.85          1        1.85
3        5      9.2           3.5          3             1        1.7
4        6      9.2           3.6          2.8           1        1.8
5        7      0             0            0             0        0
6        7.5    0             0            0             0        0
}\RGDSWtwo

\pgfplotstableread{
Location Force  GlobalSolve   GMRES    Other
1        4      7.2           6.5      0.7
2        4.5    8.2           7.4      0.8
3        5      0             0        0
4        6      0             0        0
5        7      0             0        0
6        7.5    0             0        0
}\NKS

\begin{axis}[bar shift=0pt, hide axis]
    \addplot+ table [x=Location, y=InnerSolve] {\RGDSWtwo};
    \addplot+ table [x=Location, y=CoarseSolve] {\RGDSWtwo};
    \addplot+ table [x=Location, y=GMRES] {\RGDSWtwo};
    \addplot+ table [x=Location, y=Other] {\RGDSWtwo};
\end{axis}

\begin{axis}[bar shift=11pt]
    \addplot+ table [y=InnerSolve] {\hybrid}; \addlegendentry{Inner solve}
    \addplot+ table [y=CoarseSolve] {\hybrid}; \addlegendentry{Coarse solve}
    \addplot+ table [y=GMRES] {\hybrid}; \addlegendentry{GMRES}
    \addplot+ table [y=Other] {\hybrid}; \addlegendentry{Other}
\end{axis}

\begin{axis}[bar shift=-11pt, hide axis, cycle list shift=2]
    \addplot+ table [x=Location, y=GMRES] {\NKS};
    \addplot+ table [x=Location, y=Other] {\NKS};
\end{axis}

\node[rotate=90] at (0.45,-0.6) {NKS};
\node[rotate=90] at (0.85,-0.6) {Add.};
\node[rotate=90] at (1.25,-0.6) {Hybrid};

\node[rotate=0, text width=1.5cm] (gmres) at (1.0,3.95) {GMRES its. (Newton its.)};

\node[rotate=90] (one) at (.47,2.5) {$151$ $(5)$};
\node[rotate=90] (two) at (0.89,3) {$75$ $(4)$};
\node[rotate=90] (three) at (1.28,2.45) {$44$ $(3)$};

\node[rotate=90] at (2.17,2.8) {$162$ $(5)$};
\node[rotate=90] at (2.57,3) {$76$ $(4)$};
\node[rotate=90] at (2.95,2.5) {$44$ $(3)$};

\node[rotate=0] at (3.87,0.15) {\color{red}\ding{55}};
\node[rotate=90] at (4.27,3) {$76$ $(4)$};
\node[rotate=90] at (4.68,2.4) {$44$ $(3)$};

\node[rotate=0] at (5.57,0.15) {\color{red}\ding{55}};
\node[rotate=90] at (6.0,3) {$76$ $(4)$};
\node[rotate=90] at (6.38,2.5) {$45$ $(3)$};

\node[rotate=0] at (7.28,0.15) {\color{red}\ding{55}};
\node[rotate=0] at (7.65,0.15) {\color{red}\ding{55}};
\node[rotate=90] at (8.08,2.5) {$47$ $(3)$};

\node[rotate=0] at (9.,0.15) {\color{red}\ding{55}};
\node[rotate=0] at (9.4,0.15) {\color{red}\ding{55}};
\node[rotate=0] at (9.8,0.15) {\color{red}\ding{55}};

\draw [thin] (gmres) --  (one);
\draw [thin] (gmres) --  (two);
\draw [thin] (gmres) --  (three);

\end{tikzpicture}
    \caption{Time to solution of a classical Newton-Krylov-Schwarz (NKS) solver, two-level additive and two-level hybrid nonlinear Schwarz variants for the nonlinear elasticity problem \eqref{eq:nonlinelas} solved using $576$ subdomains. The outer Newton iterations and total number of GMRES iterations summed over all outer iterations is shown above each column. A red cross indicates that the solver failed to converge.}
    \label{fig:ttselas}
\end{figure}

As a second problem that demonstrates the potential utility of nonlinear Schwarz variants we consider the beam test based on the two dimensional plane-stress Neo-Hooke model \eqref{eq:nonlinelas} introduced in section \ref{sec:intro}. Again, we use the Fritz cluster and 576 CPU cores for a METIS decomposition into 576 subdomains. We make the problem incrementally more difficult for the nonlinear solver by increasing the volume force. We use a mesh consisting of $2.6\times 10^{6}$ nodes and solve for the resulting deformation. We set the maximum number of GMRES iterations to $100$ and otherwise use the same settings for the nonlinear Schwarz variants as for the previous problem. For the Newton-Krylov-Schwarz solver we also set the relative tolerance to $10^{-4}$ and use the same settings for the GMRES solver. However, in this case we employ a (linear) two-level Schwarz preconditioner using the RGDSW coarse space.

Figure \ref{fig:ttselas} shows that the two-level nonlinear Schwarz variants are more robust as they are able to solve the problem for larger forces. The hybrid method is more robust than the additive method which coincides with the observations made previously in our MATLAB testing environment.

To summarize, we have shown some initial results from our parallel implementation of the nonlinear Schwarz variants demonstrating that it is working as expected. The weak scalability of two-level methods is superior to the classical ASPEN approach and for nonlinear elasticity problems with large deformations the nonlinear two-level Schwarz methods is more robust and just as fast as the Newton-Krylov-Schwarz approach for small test problems. Based on this implementation we plan to investigate the weak scalability of the nonlinear Schwarz variants compared to that of the classical Newton-Krylov-Schwarz method in future work.\\

\begin{acknowledgement}
The authors would like to gratefully thank  Lea Sa{\ss}mannshausen (University of Cologne) for her help using the FEDDLib software library and Sharan Nurani Ramesh (Ruhr University Bochum) for his help with the hyperelasticity problem. The authors gratefully acknowledge the scientific support and HPC resources provided by the Erlangen National High Performance Computing Center (NHR@FAU) of the Friedrich-Alexander-Universit\"at Erlangen-N\"urnberg (FAU) under the NHR project k107ce. NHR funding is provided by federal and Bavarian state authorities. NHR@FAU hardware is partially funded by the German Research Foundation (DFG) - 440719683. Finally the authors acknowledge the financial support by the German Federal Ministry of Education and Research (BMBF) in the programme SCALEXA.
\end{acknowledgement}

\bibliographystyle{spmpsci} % mathematics and physical sciences
\bibliography{NL_bib}

\end{document}